\documentclass[11pt,a4paper]{article}

\usepackage{booktabs}
\usepackage{amssymb,latexsym}
\usepackage{a4}
\usepackage[all]{xy}
\usepackage{color}

\newcommand{\U}{{\mathcal U}}
\newcommand{\UK}{E\!K}

\newcommand{\bs}{\backslash}
\newcommand{\ra}{\rangle}
\newcommand{\la}{\langle}
\newcommand{\id}{\mathrm{id}}
\newcommand{\X}{{\mathcal X}}
\newcommand{\Q}{{\mathcal Q}}
\newcommand{\R}{{\mathcal R}}

\renewcommand{\S}{{\mathcal S}}

\newcommand{\N}{{\mathbb N}}

\newcommand{\End}{{\mathrm{End}}}
\newcommand{\Hom}{{\mathrm{Hom}}}

\newcommand{\im}{{\mathrm{im}}}

\newcommand{\Gap}{{\scshape Gap}}

\newcommand{\Grig}{{\mathfrak G}}

\newcommand{\ti}{\tilde}

\newcommand{\rb}[1]{\raisebox{1.5ex}[-1.5ex]{#1}}

\usepackage{bbm}

\usepackage{rotfloat}
\floatstyle{plain}
\floatname{algorithm}{Algorithm}
\newfloat{algorithm}{tbp}{lop}

\title{{Coset enumeration for certain infinitely presented groups}}
\author{Ren\'e Hartung}
\date{May 1, 2011}

\newenvironment{proof}{\par\vskip-\lastskip\vskip\topsep
\noindent{\it Proof.}\vadjust{\nobreak}\quad
\begingroup\divide\topsep3\divide\itemsep3
\divide\partopsep3\divide\parskip3
\divide\parsep3}
{\ifvmode\penalty10000\hbox to\hsize{\hfil$\Box$}
\else\parfillskip0pt\widowpenalty10000\hfil$\Box$
\fi\par\vskip 1.5ex\endgroup}

\newtheorem{theorem}{Theorem}
\newtheorem{corollary}[theorem]{Corollary}
\newtheorem{lemma}[theorem]{Lemma}
\newtheorem{example}[theorem]{Example}

\begin{document}
\maketitle
\begin{abstract}
  We describe an algorithm that computes the index of a finitely generated
  subgroup in a finitely $L$-presented group provided that this index is
  finite. This algorithm shows that the subgroup membership problem for
  finite index subgroups in a finitely $L$-presented group is decidable.
  As an application, we consider the low-index subgroups of some
  self-similar groups including the Grigorchuk group, the twisted twin
  of the Grigorchuk group, the Grigorchuk super-group and the Hanoi
  $3$-group.\bigskip

\end{abstract}

%%%%%%%%%%%%%%%%%%%%%%%%%%%%%%%%%%%%%%%%%%%%%%%%%%%%%%%%%%%%%%%%%%%%%%%%%%%%
\section{Introduction}
Many algorithmic problems are unsolvable for finitely presented groups
in general. For instance, there is no algorithm which allows to decide
if a group given by a finite presentation is trivial~\cite{LS77}.
However, the coset enumeration process introduced by Todd \& Coxeter~\cite{TC36}
and investigated by various others, see~\cite{HR00} or the historical
notes in Chapter 5.9 of~\cite{Sims94}, computes the index
of a finitely generated subgroup in a finitely presented group provided
that this index is finite. Therefore, the Todd-Coxeter method allows one to
prove that a finitely presented group is trivial.

Coset enumeration is one of the most important tools for investigating
finitely presented groups; but, if the subgroup has infinite index,
this process will not terminate. Even if the subgroup has finite index,
there is no upper bound on the complexity of coset enumeration.  Therefore,
even proving a finitely presented group being trivial is computationally
a challenging problem~\cite{HR00,MNS79}.\smallskip

For this reason, solving algorithmic problems for infinitely presented
groups seems entirely infeasible. However, an interesting family of
recursively presented groups was recently shown to be applicable for
computer investigations. Examples of such groups arise as subgroups
of the automorphism group of a regular tree. A famous example is
the Grigorchuk group $\Grig$ which plays a prominent role in the
area of Burnside problems~\cite{Gri80}. The group $\Grig$ is finitely
generated and it admits a recursive presentation whose relations are
given recursively by the action of a finitely generated free monoid of
endomorphisms acting on finitely many relations~\cite{Lys85}. Infinite
presentations of this type are called \emph{finite $L$-presentations}
in honor of Lys\"enok's latter result for the Grigorchuk group $\Grig$;
see Section~\ref{sec:Pre} or~\cite{Bar03} for a definition.\smallskip

Finite $L$-presentations are `natural' generalizations of
finite presentations and, as the concept is quite general,
they found their application in various aspects of group theory; see, for
instance,~\cite{Bau71,KW01,OS02}. A finite $L$-presentation of a group allows
to compute its lower central series quotients~\cite{BEH08}
and the Dwyer quotients of its Schur multiplier~\cite{Har10}. The Dwyer
quotients often exhibit periodicities which yield detailed information on
the structure of the Schur multiplier in general.\smallskip

In this paper, we describe a coset enumeration process for computing
the index of a finitely generated subgroup in a finitely $L$-presented
group provided that this index is finite. In order to achieve
this method, we show in Section~\ref{sec:TCLp} that finitely many
relations are sufficient to compute an upper bound on the index
using coset enumeration for finitely presented groups. It then remains
to either prove that this upper bound is sharp or to improve the bound
otherwise. In Section~\ref{sec:Valid}, we show that the latter problem
is algorithmically decidable in general. In particular, we show that
there exists an algorithm which decides whether or not a map from the
free group over the $L$-presentations generators into a finite group
induces a homomorphism from the $L$-presented group.\smallskip

Similar to coset enumeration for finitely presented groups, our method
for finitely $L$-presented groups allows straightforward applications
including a membership test for finite index subgroups. In particular,
our method allows us to compute the number of subgroups with small index
for some self-similar groups in Section~\ref{sec:Apps}. Our explicit
computations correct the counts obtained in~\cite{BGZ03,BG02},
and hence we provide a further step towards Problem~6.1 raised
in~\cite{Gr05}.\smallskip

We have implemented our coset enumeration method and its applications
in the computer algebra system \Gap~\cite{GAP}. Computer experiments
with this implementation demonstrate that our method works reasonably
well in practice.\smallskip

In a forthcoming paper, we prove a variant of the variant of the
Reidemeister-Schreier theorem for finitely $L$-presented groups which
shows that each finite index subgroup of a finitely $L$-presented group
is finitely $L$-presented itself.

%%%%%%%%%%%%%%%%%%%%%%%%%%%%%%%%%%%%%%%%%%%%%%%%%%%%%%%%%%%%%%%%%%%%%%%%%%%%
\section{Preliminaries}\label{sec:Pre}
We briefly recall the notion of a finite $L$-presentation as
introduced in~\cite{Bar03}. For this purpose, let $F$ be a finitely
generated free group over the alphabet $\X$. Furthermore, let $\Q,\R\subset
F$ and $\Phi\subset\End(F)$ be finite subsets. Then the quadruple
$\la\X\mid\Q\mid\Phi\mid\R\ra$ is a \emph{finite $L$-presentation}. It
defines the \emph{finitely $L$-presented group}
\begin{equation} \label{eqn:LpGroup}
  G = \left\la\X~\middle|~\Q\cup\bigcup_{\sigma\in\Phi^*} 
  \R^\sigma \right\ra,
\end{equation}
where $\Phi^*$ denotes the free monoid of endomorphisms generated by
$\Phi$; that is, the closure of $\{\id\}\cup\Phi$ under composition of
endomorphisms. We will also write $G = \la\X\mid\Q\mid\Phi\mid\R\ra$
for the finitely $L$-presented group in Eq.~(\ref{eqn:LpGroup}).\medskip

Clearly, every finitely presented group $\la\X\mid\R\ra$ is finitely
$L$-presented by $\la\X\mid\emptyset\mid\emptyset\mid\R\ra$.
Therefore, finite $L$-presentations generalize the concept of finite
presentations. Other examples of finitely $L$-presented groups are
various self-similar groups or branch groups~\cite{Bar03}. For instance,
the Grigorchuk group satisfies the following 
\begin{theorem}[Lys\"enok, 1985]\label{thm:Grig}
  The Grigorchuk group $\Grig$ is finitely $L$-presented by
  \[
    \la \{a,b,c,d\} \mid \{a^2,b^2,c^2,d^2,bcd\} \mid \{\sigma\}\mid 
        \{ (ad)^4,(adacac)^4 \} \ra,
  \]
  where $\sigma$ is the endomorphism of the free group over the alphabet 
  $\{a,b,c,d\}$ induced by the mapping $a\mapsto aca$, $b\mapsto d$,
  $c\mapsto b$, and $d\mapsto c$.
\end{theorem}
\begin{proof}
  For a proof, we refer to~\cite{Lys85}.
\end{proof} 
Finite $L$-presentations are recursive presentations which are
`natural' generalizations of finite presentations. They were used by
various authors to construct groups with interesting properties; see,
for instance,~\cite{Bau71,KW01,OS02}. Furthermore, every free group
in a variety of groups that satisfies finitely many identities is finitely
$L$-presented~\cite{Bar03}; e.g., the free Burnside group $B(n,m)$
of exponent $m$ on $n$ generators is finitely $L$-presented by
\[
  \la \{a_1,\ldots,a_n\}\cup \{t\} \mid \{t\} \mid \Sigma \mid \{t^m\}\ra,
\]
where the endomorphisms $\Sigma = \{\sigma_x \mid x \in \{a_1^{\pm
1},\ldots,a_n^{\pm 1}\}\}$ are induced by the mappings
\[
  \sigma_x\colon \left\{\begin{array}{rcl@{\quad}l}
    a_i &\mapsto& a_i,& \textrm{for each }1\leq i\leq n\\
    t   &\mapsto& tx, 
  \end{array}\right. 
\]
for each $x\in\{a_1^{\pm 1},\ldots, a_n^{\pm 1}\}$.

%%%%%%%%%%%%%%%%%%%%%%%%%%%%%%%%%%%%%%%%%%%%%%%%%%%%%%%%%%%%%%%%%%%%%%%%%%%%
\section{Coset enumeration for finitely {\boldmath $L$}-presented groups}\label{sec:TCLp}
Let $G = \la \X \mid \Q \mid \Phi \mid \R \ra$ be a finitely $L$-presented
group and let $\U\leq G$ be a finitely generated subgroup with finite
index in $G$. In this section, we show that coset enumeration for
finitely presented groups yields an upper bound on the index $[G:\U]$. In
Section~\ref{sec:Valid}, it then remains to prove (or disprove) that this
upper bound is sharp.\medskip

Let $\{g_1,\ldots,g_n\}$ be a generating set for the subgroup $\U$.
We assume that the generators of $\U$
are given as words over the alphabet $\X\cup\X^{-}$. Denote the free group over
$\X$ by $F$ and let $K$ be the normal subgroup
\[
  K = \left\la \Q\cup\bigcup_{\sigma\in\Phi^*} \R^\sigma\right\ra^F.
\]
so that $G \cong F/K$ holds. Then the subgroup $E = \la g_1,\ldots,g_n\ra
\leq F$ satisfies that $\U \cong \UK/K$. Hence, we are to compute the index
$[G:\U] = [F:\UK]$.\medskip

For an element $\sigma\in\Phi^*$, we denote by $\|\sigma\|$ the usual
word-length in the generating set $\Phi$ of the free monoid $\Phi^*$. 
Define $\Phi^{(i)} = \{\sigma\in\Phi^*\mid\|\sigma\|\leq i\}$, for each
$i\in\N_0$. Then, as $\Q$, $\Phi$, and $\R$ are finite sets, the normal subgroup
\[
  K_i = \left\la \Q\cup \bigcup_{\sigma\in\Phi^{(i)}} \R^\sigma \right\ra^F
\]
is finitely generated as normal subgroup. We obtain $K =
\bigcup_{i\geq 0} K_i$ and also $\UK = \bigcup_{i\geq 0} \UK_i$. Consider the
ascending chain of subgroups
\[
  \UK_0 \leq \UK_1 \leq \ldots \leq \UK_\ell\leq \UK_{\ell+1}\leq \ldots 
  \leq \UK\leq F.
\]
Then the following lemma is straightforward.
\begin{lemma}\label{lem:FiniteGuess}
  The subgroup $\UK$ has finite index in $F$ if and only if there exists
  $\ell\in\N$ such that $\UK_\ell$ has finite index in $F$. In that 
  case, there exists $\ell'\in\N$ such that $\UK_{\ell'} = \UK$.
\end{lemma}
\begin{proof}
  Obviously, if $[F:\UK_\ell]$ is finite for some $\ell\in\N$, then 
  the subgroup $\UK$ has finite index in $F$. On the other hand,
  if $[F:\UK]$ is finite, then, as $F$ is finitely generated, the
  subgroup $\UK$ is finitely generated. Let $\{u_1,\ldots,u_n\}$ be a
  generating set of $\UK$. Since $\UK = \bigcup_{i\geq 0}
  \UK_i$ holds, there exists a positive integer $\ell\in\N$ such that
  $\{u_1,\ldots,u_n\} \subseteq \UK_{\ell}$ and
  thus $\UK_{\ell} = \UK$.
\end{proof}
Note that the index $[F:\UK_\ell]$ is the index of the subgroup
$\U$ in the finitely presented \emph{covering group}
\begin{equation}\label{eqn:FpGroup}
  G_\ell = \la\X\mid\{q,r^\sigma\mid q\in\Q,r\in\R,\sigma\in\Phi^{(\ell)}\}\ra.
\end{equation}
By Lemma~\ref{lem:FiniteGuess}, there exists a positive integer
$\ell\in\N$ so that the subgroup $\U$ has finite index in $G_\ell$. In
this case, coset enumeration for finitely presented groups 
computes the index $[G_\ell:\U]$. Although we
do not know this integer $\ell$ \emph{a priori}, we can use the
following firsthand approach to find such an integer: Starting
with $\ell=1$, we attempt to prove finiteness of $[G_\ell:\U]$
using coset enumeration for finitely presented groups. If this attempt
does not succeed within a previously defined time limit, we increase the
integer $\ell$ and the time limit. We continue this process
until eventually the index $[G_\ell:\U]$ is proved to be finite. In
theory, Lemma~\ref{lem:FiniteGuess} guarantees that this process will
terminate. Computer experiments with the implementation of our method
in \Gap\ show that this firsthand approach works reasonably well
in practice. In particular, our implementation allows to compute
the index of all subgroups considered in~\cite{Bar05,BG02,Gr05} and
Chapter~VIII of~\cite{Har00}.\medskip

Suppose that the integer $\ell\in\N$ is chosen so that $n = [G_\ell:\U]$
is finite and that the coset enumeration for finitely presented groups
has terminated and has computed a permutation representation $\varphi_\ell\colon
F\to \S_n$ for the group's action on the right cosets $\UK_\ell \bs
F$. Then the index $[G:\U] = [F:\UK]$ divides the index $[G_\ell:\U] =
[F:\UK_\ell]$, and hence $[G_\ell:\U]$ is an upper bound on $[G:\U]$. It
therefore remains to either prove that $[F:\UK] = [F:\UK_\ell]$ holds, or
to increase the integer $\ell$ otherwise. The permutation representation
$\varphi_\ell\colon F\to \S_n$ is called \emph{valid}, if $[F:\UK] =
[F:\UK_\ell]$ holds.\medskip

Clearly, a permutation representation $\varphi_\ell\colon F\to\S_n$
is valid if and only if every relation $r\in F$ of the group
presentation is contained in the kernel of $\varphi_\ell$. Therefore,
if the group $G = F/K$ were finitely presented, only finitely many
relations need to be considered to prove validity of $\varphi_\ell$.
However, for finitely $L$-presented groups, even checking validity of a
permutation representation $\varphi_\ell$ involves possibly infinitely
many relations. In Section~\ref{sec:Valid}, we prove that the latter
problem is decidable in general.

%%%%%%%%%%%%%%%%%%%%%%%%%%%%%%%%%%%%%%%%%%%%%%%%%%%%%%%%%%%%%%%%%%%%%%%%%%%%
\section{Deciding validity of a permutation representation}
\label{sec:Valid}
In this section, we describe our algorithm for deciding
whether or not a permutation representation $\varphi\colon F\to\S_n$,
as considered in Section~\ref{sec:TCLp}, is valid. This is equivalent to
checking whether a coset-table for $\U$ in $G_\ell$ obtained by the methods
of Section~\ref{sec:TCLp} defines the given subgroup $\U\leq G$.\smallskip

Let $\varphi\colon F\to \S_n$ be a permutation representation as in
Section~\ref{sec:TCLp} and let $\Phi^*$ be the free monoid generated
by a finite set $\Phi\subseteq \End(F)$. For two endomorphisms
$\sigma\in\Phi^*$ and $\delta\in\Phi^*$, we say that \emph{$\delta$
reduces to $\sigma$ with respect to $\varphi$} if there exists a
homomorphism $\pi\colon \im(\sigma\varphi) \to \im(\delta\varphi)$ such
that $\sigma\varphi\pi = \delta\varphi$. In this case, we will write
$\delta\leadsto_\varphi \sigma$. Note that $\leadsto$ is a reflexive
and transitive relation on the endomorphisms $\Phi^*$.  The following
lemma gives an equivalent definition for $\delta \leadsto_\varphi \sigma$.
\begin{lemma}
  Let $\delta,\sigma\in\End(F)$ be given. Then $\delta$ reduces to $\sigma$
  with respect to $\varphi$ if and only if $\ker(\sigma\varphi) \leq
  \ker(\delta\varphi)$ holds.
\end{lemma}
\begin{proof}
  Assume that $\delta \leadsto_\varphi \sigma$ holds. Then, by definition,
  there exists a homomorphism
  $\pi\colon \im(\sigma\varphi) \to \im(\delta\varphi)$ such that
  $\sigma\varphi\pi = \delta\varphi$. Let $g\in\ker(\sigma\varphi)$.
  Then we have that $g^{\delta\varphi} = g^{\sigma\varphi\pi}
  = (g^{\sigma\varphi})^\pi = 1$ and hence, we obtain $g\in
  \ker(\delta\varphi)$. Suppose that $\ker(\sigma\varphi)
  \leq \ker(\delta\varphi)$ holds. Then we have the isomorphisms
  $F/\ker(\sigma\varphi) \to \im(\sigma\varphi),\: g\ker(\sigma\varphi)
  \mapsto g^{\sigma\varphi}$ and $F/\ker(\delta\varphi)
  \to \im(\delta\varphi),\: g\ker(\delta\varphi) \mapsto
  g^{\delta\varphi}$. We further have the natural homomorphism
  $F/\ker(\sigma\varphi) \to F/\ker(\delta\varphi),\: g\ker(\sigma\varphi)
  \mapsto g\ker(\delta\varphi)$. This yields the existence of a homomorphism
  $\pi\colon \im(\sigma\varphi) \to \im(\delta\varphi)$ such that
  $g^{\sigma\varphi\pi} = g^{\delta\varphi}$.
\end{proof}
A finite generating set for the kernel $\ker(\sigma\varphi)$ is given
by the Schreier theorem~\cite[Proposition~3.7]{LS77} and hence, it is
straightforward to check whether or not $\delta\leadsto_\varphi \sigma$
holds. The definition $\delta \leadsto_\varphi \sigma$ also yields the 
following immediate consequence.
\begin{lemma}\label{lem:Fin}
  There is no infinite set of endomorphisms of $F$ such that for each pair
  $(\sigma,\delta)$ from this set, neither $\sigma\leadsto_\varphi \delta$
  nor $\delta\leadsto_\varphi\sigma$ hold.
\end{lemma}
\begin{proof}
  Obviously, for every endomorphism $\sigma\in\End(F)$, it holds that
  $\sigma \leadsto_\varphi \sigma$. By the universal property of the free
  group $F$, a homomorphism $\sigma\varphi\colon F\to \S_n$ is uniquely
  defined by the images $x_1^{\sigma\varphi},\ldots,x_n^{\sigma\varphi}$
  of the elements $x_1,\ldots,x_n$ of a basis of $F$. Since
  $\im(\varphi)$ is a finite group, there are only finitely many
  homomorphisms $F\to\im(\varphi)$ and therefore $\Hom(F,\im(\varphi))$
  is finite. Hence, an infinite set of endomorphisms will contain
  endomorphisms $\sigma$ and $\delta$ with $x_i^{\sigma\varphi}
  = x_i^{\delta\varphi}$, for each $1\leq i\leq n$. In this case,
  $\sigma\leadsto_\varphi \delta$ obviously holds.
\end{proof}
An element $\sigma\in\Phi^*$ is called a \emph{$\Phi$-descendant
of $\delta\in\Phi^*$} if there exists $\psi\in\Phi$ such that $\sigma
= \psi\delta$. Thereby, the free monoid $\Phi^*$ obtains the
structure of a $|\Phi|$-regular rooted tree with its root being the
identity map $\id\colon F\to F$. We can further endow the monoid $\Phi^*$
with a length-plus-(from the right)-lexicographic ordering $\prec$ by
choosing an arbitrary ordering on the finite set $\Phi$. More precisely,
we define $\sigma\prec\delta$ if $\|\sigma\|<\|\delta\|$ or, otherwise, if
$\sigma = \sigma_1\cdots \sigma_n$ and $\delta = \delta_1\cdots\delta_n$,
with $\sigma_i,\delta_i\in\Phi$, and there exists a positive integer
$1\leq k\leq n$ such that $\sigma_i = \delta_i$, for $k<i\leq n$, and
$\sigma_k\prec\delta_k$. Since $\Phi$ is finite, the obtained ordering
$\prec$ is a well-ordering on the monoid $\Phi^*$, see~\cite{Sims94}, and
therefore there is no infinite $\prec$-descending series of endomorphisms
in $\Phi^*$.\smallskip

Our algorithm that decides validity of a permutation
representation $\varphi\colon F\to\S_n$ is displayed in
Algorithm~\ref{alg:IsValidPermRep} below.
\begin{algorithm}[htp]
  \caption{Deciding validity of a permutation representation}
  \label{alg:IsValidPermRep}
  \begin{center}
  {\begin{minipage}{10cm}
    \begin{tabbing}
      {\scshape IsValidPermRep}($\X$, $\Q$, $\Phi$, $\R$, $\U$, $\varphi$)\\
     \quad\= Initialize $V:=\{\id\colon F\to F\}$ and $S:=\Phi$.\\
      \> Choose an ordering on $\Phi = \{\phi_1,\ldots,\phi_n\}$ with $\phi_i \prec \phi_{i+1}$.\\
      \>{\bf while} $S \neq \emptyset$ {\bf do} \\
      \>\quad\= Remove the first entry $\delta$ from $S$.\\
      \>\> {\bf if} $\left( \exists\,r\in\R\colon\: r^{\delta} \not\in \ker\varphi \right)$ {\bf then} {\bf return}( \verb#false# )\\
      \>\>{\bf if not} $\left( \exists\,\sigma\in V\colon\: \delta\leadsto_\varphi\sigma \right)$ {\bf then}\\
      \>\>\quad\= Append $\phi_1\delta,\ldots,\phi_n\delta$ to $S$.\\
      \>\>\> Add $\delta$ to $V$. \\
      \>{\bf return}( \verb#true# )
    \end{tabbing}
  \end{minipage}}
  \end{center}
\end{algorithm}
We need to prove the following 
\begin{theorem}\label{thm:Validity}
  The algorithm {\scshape IsValidPermRep} returns \verb|true| if and only
  if the permutation representation $\varphi\colon F\to\S_n$ is valid. 
\end{theorem}
\begin{proof} 
  The ordering $\prec$ on $\Phi$ can be extended to an ordering on $\Phi^*$ 
  as described above. By construction, the stack $S$ is ordered with respect to $\prec$. 
  Since $F$ is finitely generated, the set of homomorphisms
  $\Hom(F,\S_n)$ is finite. Thus, in particular, the set $\{\delta\varphi
  \mid \delta\in V\} \subseteq \Hom(F,\S_n)$ is finite and therefore
  the algorithm {\scshape IsValidPermRep} can add only finitely many endomorphisms to the
  set $V$. Thus, for every $\Phi$-descendant $\delta$ in the stack $S$,
  there will eventually exist an element $\sigma\in V$ such that $\delta
  \leadsto_\varphi \sigma$. Therefore, the algorithm {\scshape
  IsValidPermRep} is guaranteed to terminate and it returns either
  \verb|true| or \verb|false|. Clearly, if the algorithm returned
  \verb#false#, then it found a relation $r^\delta$ which yields a
  coincidence, and hence the permutation representation $\varphi\colon
  F\to\S_n$ is not valid.

  Suppose that the algorithm returned \verb#true#. By the constructions
  of Section~\ref{sec:TCLp}, the fixed relations in $\Q$ of the 
  $L$-presentation $\la\X\mid\Q\mid\Phi\mid\R\ra$ are
  already contained in the kernel of the permutation representation
  $\varphi$. Therefore, it suffices to prove that every relation of
  the form $r^{\sigma_1}$, with $r\in\R$ and $\sigma_1\in\Phi^*$, is
  contained in the kernel of $\varphi$. By construction, there exists
  $\delta\in V$ maximal subject to the existence of $w\in\Phi^*$
  such that $\sigma_1 = w\delta$. If $\|w\| = 0$, then $\sigma_1 =
  \delta$ is contained in $V$ and therefore $r^\delta\in\ker\varphi$,
  as the algorithm did not return \verb#false#. Otherwise,
  there exist $\psi\in\Phi$ and $v\in\Phi^*$ such that $w\delta =
  v\psi\delta$. Since $\psi\delta\not\in V$, there exists
  an element $\varepsilon\in V$ with $\varepsilon\prec\psi\delta$,
  by construction, such that $\psi\delta$ reduces to $\varepsilon$
  with respect to $\varphi$. Thus, by definition, there exists
  a homomorphism $\pi\colon \im(\varepsilon\varphi)
  \to \im(\psi\delta\varphi)$ such that $\psi\delta\varphi =
  \varepsilon\varphi\pi$. In particular, we obtain that
  $r^{\sigma_1\varphi} = r^{w\delta\varphi} = r^{v\psi\delta\varphi}
  = r^{v\varepsilon\varphi\pi}$. As $\pi$
  is a homomorphism, it suffices to prove that $r^{v\varepsilon}\in
  \ker\varphi$. Note that, since $\varepsilon\prec \psi\delta$, we
  have that $v\varepsilon\prec v\psi\delta = \sigma_1$. Continuing this
  rewriting process with the element $\sigma_2 = v\varepsilon$ yields a
  descending sequence $\sigma_1 \succ\sigma_2\succ\ldots$ in the monoid
  $\Phi^*$. As the ordering $\prec$ is a well-ordering, this process
  terminates with an element $\sigma_n\in V$. Since the algorithm did
  not return \verb#false#, we have that $r^{\sigma_n}\in\ker\varphi$
  which proves the assertion.
\end{proof}
Note that, if the algorithm {\scshape IsValidPermRep} found a coincidence,
this can be used to update the coset-table, and thus another application of
coset enumeration for finitely presented groups can be avoided. Moreover,
the Algorithm~\ref{alg:IsValidPermRep} yields the following
\begin{theorem}
  Let $G$ be finitely $L$-presented by $\la \X\mid \Q\mid\Phi\mid\R\ra$
  and denote the free group over $\X$ by $F$. There exists an
  algorithm which decides whether or not a homomorphisms $\varphi\colon
  F \to \S_n$ induces a homomorphism $G \to \S_n$.
\end{theorem}
If $\Phi^* = \{\sigma\}^*$ is generated by a single element
$\sigma\in\End(F)$, then there will exist positive integers $0\leq i< j$ such
that $\sigma^j \leadsto_\varphi \sigma^i$. In this case, the algorithm
{\scshape IsValidPermRep} simplifies to the following
\begin{corollary}\label{cor:CondCyc}
  Let $0\leq i<j$ be positive integers such that $\sigma^j\leadsto_\varphi 
  \sigma^i$. Then we have $[F:\UK_\ell] = [F:\UK]$ if and only if
  \begin{equation} \label{eqn:FinRels}
    \{q,r^{\sigma^k}\mid  q\in\Q,\:r\in\R,\:0\leq k< j\} 
    \subseteq \ker\varphi.
  \end{equation}
\end{corollary}
We consider the following 
\begin{example}\label{ex:Bas}
  Let $G$ denote the Basilica Group~\cite{GZ02b}. Then
  $G$ is finitely $L$-presented by $\la \{a,b\}
  \mid \emptyset \mid \{\sigma\} \mid \{[a,a^b]\} \ra$, where
  $\sigma$ is induced by the mapping $a\mapsto b^2$ and $b\mapsto a$;
  see~\cite{BV05}. We consider the subgroup $\U = \la  a^3, b, aba \ra$.
  A coset enumeration for finitely presented groups yields that the subgroup
  $\U$ has index $3$ in the finitely presented covering group
  \[
    G_0 = \la \{a,b\} \mid \{[a,a^b]\} \ra.
  \]
  Furthermore, we obtain the permutation representation $\varphi\colon
  F\to\S_3$ for the group's action on the cosets $\UK_0\backslash F$. This 
  permutation representation is induced by the mapping
  \[
    a \mapsto (1,2,3)\quad\textrm{and}\quad
    b \mapsto (2,3).
  \]
  We now obtain the images 
  \[
    \begin{array}{rcl@{\qquad}rcl}
      a^{\sigma\varphi} &=& (\:),& 
      b^{\sigma\varphi} &=& (1,2,3),\\
      a^{\sigma^2\varphi}&=& (1,3,2),&
      b^{\sigma^2\varphi}&=& (\:),\\
      a^{\sigma^3\varphi}&=& (\:),& 
      b^{\sigma^3\varphi}&=& (1,3,2).
    \end{array}
  \]
  Clearly, the mapping $a^{\sigma\varphi} \mapsto a ^ {\sigma^3\varphi}$
  and $b^{\sigma\varphi} \mapsto b ^ {\sigma^3\varphi}$
  induces a homomorphism $\pi\colon \im(\sigma\varphi) \to
  \im(\sigma^3\varphi)$, and hence we have $\sigma^3 \leadsto_\varphi
  \sigma$. By Corollary~\ref{cor:CondCyc}, it therefore suffices to
  prove that
  \[
    ([a,a^b])^ {\varphi} = (\:),\quad
    ([a,a^b])^ {\sigma\varphi} = (\:),\quad\textrm{and}\quad
    ([a,a^b])^ {\sigma^2\varphi} = (\:)
  \]
  hold. This yields that $[G:\U] = 3$.\smallskip
\end{example}

%%%%%%%%%%%%%%%%%%%%%%%%%%%%%%%%%%%%%%%%%%%%%%%%%%%%%%%%%%%%%%%%%%%%%%%%%%%%
\section{Further applications}\label{sec:Apps}
The permutation representation $\varphi\colon F\to\S_n$ for a finite index
subgroup $\UK/K \leq F/K$ yields various algorithmic applications. For
instance, an element $w\in F$ is contained in the given subgroup $\UK$
if and only if it stabilizes the trivial coset $\UK\,1$. This can be
easily be checked using the permutation representation $\varphi$. In
particular, we obtain
\begin{theorem}
  The subgroup membership problem for finite index subgroups in a finitely 
  $L$-presented group is decidable.
\end{theorem}
Moreover, having computed permutation representations $\varphi_1$ and
$\varphi_2$ for two finite index subgroups $U$ and $V$ of a finitely
$L$-presented group, one can compute a generating set for the intersection
$U\cap V$. Thus, in particular, our method allows one to compute the core
of a finite index subgroup. For example, the core of the subgroup $\U$
in Example~\ref{ex:Bas} is given by
\[
  H = \la b^2, a^3, a^2ba^{-1}b^{-1}, abab^{-1}, ab^2a^{-1},
          ba^2b^{-1}a^{-1}, baba^{-2} \ra.
\]
Since $H$ has finite index in $G$, our method allows to compute a
permutation representation for the core $H$ and we obtain $G / H \cong
\S_3$.

\subsection{Low-index subgroups of finitely $L$-presented groups}
The coset enumeration process for finitely presented groups was
used in~\cite{DSch74} to describe a low-index subgroup algorithm that
computes all subgroups of a finitely presented group up to a given
index. This algorithm also yields a method for computing all subgroups
with small index in a finitely $L$-presented group. In this section,
we will describe this method for finitely $L$-presented groups and we
use this algorithm to investigate some self-similar groups. In
particular, our implementation in the computer algebra system \Gap\
allows us to determine the number of subgroups with index at most $64$
in the Grigorchuk group.\medskip

Let $G = F/K$ be a finitely $L$-presented group and let $n$ be a
positive integer. Using the low-index subgroup algorithm for finitely
presented groups~\cite{DSch74}, see also Chapter~5.6 of~\cite{Sims94},
we obtain the list of subgroups with index at most $n$ in the finitely
presented covering group $G_\ell = F/K_\ell$. Since the covering group
$G_\ell$ naturally maps onto $G$, every subgroup $\UK_\ell/K_\ell$ with
index $n$ in $G_\ell$ maps to a subgroup of the finitely $L$-presented
group $G$. The index of this image $\UK/K$ in $G$ divides the index $n =
[F:\UK_\ell]$.  On the other hand, every subgroup $\UK/K$ with index $n$
in the finitely $L$-presented group $F/K$ has a full preimage $\UK/K_\ell$
in the finitely presented group $G_\ell$ with index $n$. Thus the list of
subgroups with index at most $n$ in a finitely $L$-presented group $G$ can
be obtained from the list of subgroups of a finitely presented covering
group $G_\ell$ by removing duplicate images. Our solution to the subgroup
membership problem can be used to remove duplicate images in $G$.\medskip

As an application, we consider some interesting self-similar groups
and we determine the number of subgroups with small index. We first
consider the Grigorchuk group $\Grig$: its lattice of normal subgroups is
well-understood~\cite{Bar05} while its lattice of subgroups with finite
index is widely unknown~\cite{Gr05}. It is well known~\cite{Gr05} that
the Grigorchuk group has seven subgroups of index two. In~\cite{Per00}, it
was shown that the subgroups of index two are the only maximal subgroups
of $\Grig$. Our low-index subgroup algorithm allows us to determine
the number of subgroups with index at most $64$ in the group $\Grig$
and thereby, it corrects the counts in Section~7.4 of~\cite{BGZ03} and in Section~4.1 of~\cite{BG02}. The following
list summarizes the number of subgroups ($\leq$) and the number of
normal subgroups ($\unlhd$) among them:
\[
  \begin{array}{cccccccc}
    \toprule
    {\rm index} & 1 & 2 & 4 & 8 & 16 & 32 & 64\\ 
    \midrule
     \leq       & 1 & 7 & 31 & 183 & 1827 & 22931 & 378403 \\
     \unlhd     & 1 & 7 & 7  & 7   & 5    & 3     & 3 \\
    \bottomrule
  \end{array}
\]
The Grigorchuk super-group $\ti\Grig$ was introduced in~\cite{BG02}. It
contains the Grigorchuk group as an infinite index subgroup. Little
is known about its subgroup lattice. The twisted twin $\bar\Grig$ of the 
Grigorchuk group was introduced in~\cite{BS10}. Similarly, little is known about
the subgroup lattice of the twisted twin $\bar\Grig$. 
Our low-index subgroup algorithm allows us to determine
the number of subgroups with index at most $16$ in both groups. Their
subgroup counts are:
\[
  \begin{array}{crrrr}
   \toprule
                  & \multicolumn{2}{c}{\ti\Grig} & \multicolumn{2}{c}{\bar\Grig} \\
   \rb{\rm index} & \multicolumn{1}{c}{  \leq   } & \multicolumn{1}{c}{ \unlhd   } & \multicolumn{1}{c}{ \leq    } & \multicolumn{1}{c}{  \unlhd  } \\
    \midrule
    1 &      1 & 1  &     1 &  1 \\
    2 &    15  &15  &    15 & 15 \\
    4 &   147  &35  &   147 & 35 \\
    8 &  2163  &43  &  1963 & 43 \\
   16 & 52403  &55  & 46723 & 47 \\
   \bottomrule
  \end{array}
\]
As both groups are $2$-groups, the only maximal subgroups with finite index are the subgroups with index two; though the 
question of dertermining all maximal subgroups of $\ti\Grig$ and $\bar\Grig$ has not been addressed in this paper.\smallskip

Finally, we consider the Basilica group and the 
Hanoi-$3$ group~\cite{GS06} with its $L$-presentation from~\cite{BSZ10}. 
The following list also includes the number of maximal subgroups (max):
\[
    \begin{array}{crrrrrr}
    \toprule
    &\multicolumn{3}{c}{{\rm Hanoi-}3} & \multicolumn{3}{c}{\rm Basilica}\\ 
    \rb{\rm index} & \leq & \unlhd & {\rm max} & \leq & \unlhd & {\rm max}  \\
    \midrule
    1 &    1& 1  & 1   &    1 &  1 & 1\\
    2 &    7& 7  & 7   &    3 &  3 & 3\\
    3 &   12& 0  & 12  &    7 &  4 & 7\\
    4 &   59& 7  & 4   &   19 &  7 & 0\\
    5 &   15& 0  & 15  &   11 &  6 &11\\
    6 &  136& 4  & 0   &   39 & 13 & 0\\
    7 &   21& 0  & 21  &   15 &  8 &15\\
    8 &  335& 13 & 0   &  163 & 19 & 0\\
    9 &  225& 0  & 0   &  115 & 13 & 9\\
   10 &  153& 3  & 0   &   83 & 19 & 0\\
   11 &   33& 0  & 33  &   23 & 12 &23\\
   12 & 2872& 12 & 0   &  355 & 31 & 0\\
   13 &   39& 0  & 39  &   27 & 14 &27\\
   14 &  297& 3  & 0   &  115 & 25 & 0\\
   15 &  450& 0  & 0   &   77 & 24 & 0\\
   16 & 1855& 13 & 0   & 1843 & 47 & 0\\
   17 &   51& 0  & 51  &   35 & 18 &35\\
   18 & 5001& 3  & 0   & 1047 & 44 & 0\\
   19 &   57& 0  & 57  &   39 & 20 &39\\
   20 & 1189& 9  & 0   &  939 & 45 & 0\\
   21 &  756& 0  & 0   &  105 & 32 & 0\\
   22 &  531& 3  & 0   &  223 & 37 & 0\\
   23 &   69& 0  &69   &   47 & 24 &47\\
   24 &52220& 23 & 0   & 4723 & 87 & 0\\
   25 &  225& 0  &75   &  411 & 31 &25\\
   26 &  783& 3  & 0   &  315 & 43 & 0\\
   27 & 5616& 0  &27   &  736 & 49 & 0\\
   28 & 2301& 9  & 0   &      &    & \\
   29 &   87& 0  &87   &      &    & \\
   30 &15462& 3  & 0   &      &    & \\
   31 &   93& 0  &93   &      &    & \\
   32 & 9119& 25 & 0   &      &    & \\
   \bottomrule
   \end{array}
\]
The largest abelian quotient $H/H'$ of the Hanoi-$3$ group $H$
is $2$-elementary abelian of rank $3$. Thus, by the Feit-Thompson theorem~\cite{FT63},
there are no normal subgroups with odd index in the Hanoi-$3$ group.

\subsection*{Acknowledgments}
I am grateful to Laurent Bartholdi for valuable comments and suggestions.

\def\cprime{$'$}

\noindent Ren\'e Hartung,
{\scshape Mathematisches Institut},
{\scshape Georg-August Universit\"at zu G\"ottingen},
{\scshape Bunsenstra\ss e 3--5},
{\scshape 37073 G\"ottingen}
{\scshape Germany}\\[1ex]
{\it Email:} \qquad \verb|rhartung@uni-math.gwdg.de|\\[2.ex]
June 2010 (revised May 2011).
\end{document}